\documentclass{amsart}

\usepackage{amssymb,amsbsy,amsthm,amsmath,graphicx,epsfig,color}

\newtheorem{thm}{Theorem}[section]
\newtheorem*{thm*}{Theorem}

\newtheorem{conj}[thm]{Conjecture}

\theoremstyle{remark}

\newtheorem*{remark}{Remark}

\newcommand{\T}{\mathbb{T}}

\newcommand{\N}{\mathbb{N}}

\newcommand{\C}{\mathbb{C}}

\newcommand{\Z}{\mathbf{Z}}

\newcommand{\veps}{\varepsilon}

\newcommand{\aveN}{\frac{1}{N}\sum_{n=1}^N}
\newcommand{\limaveN}{\lim_{N\to\infty} \aveN}

\newcommand{\ol}[1]{\overline{#1}}

\newcommand{\bcdot}{{\bullet}}

\title[A polynomial version of Sarnak's conjecture]{A polynomial version of Sarnak's conjecture
} 

\author{Tanja Eisner}
\address{Institute of Mathematics, University of Leipzig,
P.O. Box 100 920, 04009 Leipzig, Germany}
\email{eisner@math.uni-leipzig.de}



\begin{document}

\maketitle

\begin{abstract}
Motivated by the variations of Sarnak's conjecture due to El Abdalaoui, Kulaga-Przymus, Lema\'nczyk, de la Rue and by the observation that the M\"obius function is a good weight (with limit zero) for the polynomial pointwise ergodic theorem in $L^q$, $q>1$, we introduce polynomial versions of the Sarnak conjecture for minimal systems.



\end{abstract}

\section{Introduction}

A recent conjecture of Sarnak, see \cite{sarnak-lectures}, states that the M\"obius function
$$
   \mu (n):= \begin{cases} (-1)^k, \quad & n \text{ is the product of } k \text{ distinct primes}, \\ 0, \quad & n \text{ is not square free} \end{cases}
$$
is asymptotically orthogonal  to (or
disjoint from) every deterministic sequence $(a_n)\subset \C$, i.e., the sequence $(\mu(n))$ satisfies
\begin{equation}\label{eq:sarnak}
  \limaveN \mu(n) a_n = 0
\end{equation}
for every $(a_n)$
obtained as
$a_n=f(T^n x)$
for some compact metric space $X$, $T:X\to X$ continuous with topological entropy $0$, $x\in X$ and $f\in C(X)$.
In the measure-theoretic setting the conjecture is true, i.e., for every standard probability space $(X,\nu)$, every $\nu$-preserving $T:X\to X$ and every $f\in L^1(X,\nu)$,
$$
\limaveN \mu(n) f(T^n x) =  0
$$
holds for \emph{a.e.}~$x\in X$. This means that the sequence $(\mu(n))$ is a good weight for the pointwise ergodic theorem and was noted by Sarnak in \cite{sarnak-lectures}, see El Abdalaoui, Kulaga-Przymus, Lema\'nczyk, de la Rue \cite[Section 3]{AKLR}. See also Cuny, Weber \cite[Section 2]{CW} for related results for the M\"obius and the Liouville function.

 Sarnak's conjecture follows from the Chowla conjecture (see \cite{sarnak-lectures}, \cite{AKLR} or Tao \cite{TT-blog}) and has been proven for certain systems, see e.g.~\cite{TT-blog}, Bourgain, Sarnak, Ziegler \cite{BSZ} and \cite{AKLR} for examples and references.
We mention here that the prime number theorem (in arithmetic progressions) implies the conjecture for constant (resp., periodic) sequences, the classical Davenport estimate
\begin{equation}\label{eq:dave}
\sup_{\lambda\in \T} \left| \aveN \mu(n) \lambda^n \right| = O_A(\log^{-A}N) \quad \forall A>0,
\end{equation}
$\T$ being the unit circle, covers quantitatively the case of
sequences of the form $(\lambda^n)$, $\lambda\in \T$, whereas general nilsequences are covered quantitatively by a generalization of the Davenport estimate to nilsequences due to Green, Tao \cite{GT12}, see Theorem \ref{thm:GT-nil} below (see also Green, Tao \cite{GT06} for the quadratic case) for connected groups and qualitatively by Ziegler \cite{Z} for possibly non-connected groups.

We also mention the so-called K\'atai-Bourgain-Sarnak-Ziegler orthogonality criterion as a sufficient condition for systems to satisfy the Sarnak conjecture, see Bourgain, Sarnak, Ziegler \cite{BSZ} and Tao \cite{TT-blog}.

Several variations of Sarnak's conjecture can be found in El Abdalaoui, Kulaga-Przymus, Lema\'nczyk, de la Rue \cite{AKLR}. Another version of Sarnak's conjecture was formulated by Tao \cite{TT-blog}: Equality (\ref{eq:sarnak}) holds for every sequence $(a_n)\subset\C$ with topological entropy zero. We refer to \cite{TT-blog} for the notion of topological entropy of a sequence and mention here that for sequences of the form $a_n=f(T^nx)$ coming  from a minimal topological dynamical system it is equivalent to the usual notion of topological entropy.

The aim of this note is to introduce several polynomial versions of Sarnak's conjecture for minimal systems.

\textbf{Acknowledgements.}
The construction in Remark  \ref{rem:frantz} is due to Nikos Frantzikinakis. The author is very grateful to him and thanks  B\'alint Farkas and Mariusz Lema\'nczyk for several helpful discussions.

\section{Polynomial orbits}

The starting point for our considerations is the following result, see Green, Tao \cite{GT10} for the terminology.
(Here $G$ is assumed to be connected and simply connected.)
\begin{thm}[Nilsequence version of Davenport's estimate, see Green, Tao \cite{GT12}]\label{thm:GT-nil}
Let $G/\Gamma$ be a nilmanifold of dimension $m\geq 1$ with some filtration $G_\bcdot$ of degree $d\geq 1$ and the metric induced by some $Q$-rational Mal'cev basis with $Q\geq 2$ adapted to the filtration $G_\bcdot$. Let further $g\in \mathrm{poly}(\Z,G_\bcdot)$ be a polynomial sequence and $F:G/\Gamma\to [-1,1]$ be Lipschitz. Then
$$
\aveN \mu(n) F(g(n)\Gamma) = O_{m,d,A}(1) Q^{O_{m,d,A}(1)} (1+\|F\|_{\mathrm{Lip}}) \log^{-A} N
$$
holds for every $A>0$.
\end{thm}
For an integer polynomial $p$,  applying this result to $G:=\T$, $\Gamma:=\{1\}$, the filtration $G_\bcdot$ of degree $\deg p+1$ given by $\T=\T=\ldots=\T\supset\{1\}$, the polynomial sequence $g(n):= \lambda^{p(n)}$ for $\lambda\in \T$ and $F_1,F_2$ with $F_1(\lambda):=\mathrm{Re} \lambda$ and $F_2(\lambda):=\mathrm{Im} \lambda$ yields in particular
$$
\sup_{\lambda\in\T} \left|\aveN \mu(n) \lambda^{p(n)}\right| = O_{A}(\log^{-A} N)
$$
for every $A>0$, being a polynomial analogue of (\ref{eq:dave}).

Combining this estimate and adapting the proof that $(\mu(n))$ is a
good weight for the classical pointwise ergodic theorem, see \cite[Proof of Prop. 3.1]{AKLR}, to the polynomial case
we obtain the following.
\begin{thm}\label{thm:good-weight}
Let $(X,\mu,T)$ be a standard invertible measure-theoretic dynamical system, $f\in L^q(X,\mu)$ for some $q>1$, and let $p$ be an integer polynomial. Then
$$
\limaveN \mu(n) f(T^{p(n)} x)  = 0
$$
holds for a.e. $x$. In particular, the sequence $(\mu(n))$ is a good weight for the polynomial pointwise ergodic theorem in $L^q$ for every $q>1$.
\end{thm}
\begin{remark}
\begin{itemize}
\item[1)] Using invertible extensions, the assertion also holds for non-invertible systems and
polynomials $p:\N\to\N_{0}$.
%
\item[2)] By Chu \cite{chu}, Theorem \ref{thm:GT-nil} and the argument for possibly non-connected groups by Ziegler \cite{Z}, $(\mu(n))$ is also a good weight (with limit $0$) for the multiple polynomial mean ergodic theorem. We refer to Frantzikinakis, Host \cite[Thm. 1.1]{FH} for a different argument for commuting transformations.
\end{itemize}
\end{remark}

Thus, the following variation of Sarnak's conjecture arises naturally.
\begin{conj}[Polynomial Sarnak's conjecture for minimal systems I]
Let $(X,T)$ be a minimal topological dynamical system with entropy zero. Then
\begin{equation}\label{eq:sarnak-pol}
\limaveN \mu(n) f(T^{p(n)} x)  = 0
\end{equation}
holds for every $f\in C(X)$, every polynomial $p:\N \to \N_0$ and \emph{every} $x\in X$.
\end{conj}

We present a variant of the above conjecture.
By systems with polynomial entropy zero we mean systems with sequential topological entropy equal to zero along all polynomial sequences $(p(n))$, $p$ being a polynomial mapping $\N$ to $\N_0$. For the definition and basic properties of topological entropy along sequences  we refer to Goodman \cite{G}.

\begin{conj}[Polynomial Sarnak's conjecture for minimal systems II]
Let $(X,T)$ be a minimal topological dynamical system with polynomial entropy zero. Then (\ref{eq:sarnak-pol})
holds for every $f\in C(X)$, every polynomial $p:\N \to \N_0$ and every $x\in X$.
\end{conj}


\begin{remark}\label{rem:frantz}
As observed by Nikos Frantzikinakis and Mariusz Lema\'nczyk, the conjectures I-II are false without the minimality assumption. A counterexample consists of the polynomial $p(n):=n^2$, the subshift of $\{-1,0,1\}^{\N_0}$ generated by $a:=(a_n)$ with
$$
a_n:=\begin{cases}
\mu(k), \quad &\ n=k^2 \text{ for some } k\in \N,\\
0 \quad &\text{ otherwise}
\end{cases}
$$
and $f(x):=x_0$. In this case
$$
\limaveN \mu(n) f(T^{n^2} a) =\limaveN \mu(n)^2=\frac{6}{\pi^2}
$$
holds.
On the other hand, for the standard metric and every $k\in\N$ we have that for $\veps:=2^{-k}$, the maximal number $N(n,\veps)$ of $(n,\veps)$-distinguishable orbits is bounded by the position of the beginning of the first block of zeros in $a$ of length $n+k$  plus the number of unit vectors in $\C^{n+k}$. Therefore  $N(n,\veps)$ is bounded by  a quadratic polynomial in $n$, implying
zero topological entropy of the system along every polynomial.
\end{remark}

Finally we introduce one more natural polynomial variation of Sarnak's conjecture in the spirit of Conjecture 2 in Tao \cite{TT-blog} (see there for the abstract notion of topological entropy for complex sequences).
 \begin{conj}[Polynomial Sarnak's conjecture for minimal systems III]
Let $p:\N \to \N_0$ be a polynomial and let $(X,T)$ be a minimal topological dynamical system.
Then (\ref{eq:sarnak-pol}) holds for every $f\in C(X)$ and every $x\in X$ such that the sequence $(f(T^{p(n)}x))$ has topological entropy zero.
\end{conj}

Since polynomial sequences are highly structured, the above conjectures fit well with the general intuition that the M\"obius function is asymptotically ortho\-gonal to structured sequences. Furthermore, considering averages of polynomial orbits of dynamical systems has long history in ergodic theory and combinatorics, see e.g.~Furstenberg \cite{F}, Bourgain \cite{B}, Bergelson, Leibman \cite{BL-Sz} and Green, Tao \cite{GT10}. As the celebrated polynomial ergodic theorem of Bourgain \cite{B} illustrates, showing pointwise convergence along polynomial subsequences can be much more challenging than showing pointwise convergence along the natural numbers.

\section{Remarks}

Clearly, Conjecture I implies Conjecture II and Conjecture III follows from the abstract Sarnak conjecture given in \cite{TT-blog}. However, the converse implications and the connection to Sarnak's and Chowla's conjecture are not clear to us.

Since for a nilsequence $(f(a^nx))$, the polynomial nilsequence $(f(a^{p(n)}x))$ is again a nilsequence, see Chu \cite[Prop. 2.1]{chu}, Leibman \cite{L} and Green, Tao \cite[Prop. C.2]{GT10}, the polynomial Sarnak conjectures are satisfied for transitive (equivalently, minimal) nilsystems by Ziegler's argument  \cite{Z}. For connected groups one can apply Theorem  \ref{thm:GT-nil} to the sequence  $(f(a^{p(n)}x))$ directly, using equidistribution of polynomial sequences, see Leibman \cite{L}, for the standard limiting argument to reach all continuous functions. Note that transitive nilsystems have zero polynomial entropy, see Host, Kra, Maass \cite{HKM}.

Note finally that the K\'atai-Bourgain-Sarnak-Ziegler orthogonality criterion, see Bourgain, Sarnak, Ziegler \cite{BSZ} and Tao \cite{TT-blog}, implies that the property
$$
\limaveN f(T^{p(nq_1)}x) \ol{f(T^{p(nq_2)}x)} =0 \quad \text{ for all distinct primes } q_1,q_2
$$
is sufficient for (\ref{eq:sarnak-pol}).


\begin{thebibliography}{10}




\bibitem{AKLR} H.~El Abdalaoui, J.~Kulaga-Przymus, M.~Lema\'nczyk, T.~de la Rue,
\emph{The Chowla and the Sarnak conjectures from ergodic theory point of view}, preprint, available at http://arxiv.org/abs/1410.1673.

\bibitem{BL-Sz} V.~Bergelson, A.~Leibman, \emph{Polynomial extensions of van der Waerden's and Szemeredi's theorems}, Journal of AMS \textbf{9} (1996), 725--753.

\bibitem{BL} V.~Bergelson, A.~Leibman, \emph{Distribution of values of bounded generalized polynomials}, Acta Mathematica \textbf{198} (2007), 155--230.

\bibitem{B} J.~Bourgain, \emph{Pointwise ergodic theorems for arithmetic sets}, Inst. Hautes \'Etudes Sci. Publ. Math. \textbf{69} (1989), 5--45. With an appendix by the author, Harry Furstenberg, Yitzhak Katznelson and Donald S. Ornstein.

\bibitem{BSZ} J.~Bourgain, P.~Sarnak, T.~Ziegler, \emph{Distjointness of M\"obius from horocycle flows}, preprint, http://arxiv.org/abs/1110.0992.

\bibitem{chu}
Q.~Chu, \emph{Convergence of weighted polynomial multiple ergodic averages}, Proc. Amer. Math. Soc. \textbf{137} (2009), 1363--1369.

\bibitem{CW} C.~Cuny, M.~Weber, \emph{Ergodic theorems with arithmetic weights},
preprint, available at http://arxiv.org/abs/1412.7640.


\bibitem{D} H.~Davenport, \emph{On some infinite series involving arithmetical functions. II}, Quart. J. Math. Oxf.
\textbf{8} (1937), 313--320.

\bibitem{FH} N. Frantzikinakis, B. Host, \emph{Multiple ergodic theorems for arithmetic sets}, preprint, available at http://arxiv.org/abs/1503.07154.




\bibitem{F} H.~Furstenberg, \emph{Ergodic behavior of diagonal measures and a theorem of Szemer\'edi on arithmetic progressions}, J. Analyse Math. \textbf{31} (1977), 204--256.

\bibitem{G} T.~N.~T.~Goodman, \emph{Topological sequence entropy}, Proc. London Math. Soc.   \textbf{29} (1974), 331--350.

\bibitem{GT06} B.~Green, T.~Tao, \emph{Quadratic uniformity of the M\"obius function}, Annales de l'institut Fourier \textbf{58} (2008), 1863--1935.

\bibitem{GT10} B.~Green, T.~Tao, \emph{The quantitative behaviour of polynomial orbits on nilmanifolds}, Ann. of Math. (2)  \textbf{175} (2012), 465--540.

\bibitem{GT12} B.~Green, T.~Tao, \emph{The M\"obius function is strongly orthogonal to nilsequences},
Ann. of Math. (2) \textbf{175} (2012), 541--566.

\bibitem{HKM} B.~Host, B.~Kra, A.~Maass, \emph{Complexity of nilsystems and systems lacking nilfactors}, J. Anal. Math. \textbf{124} (2014), 261--295.

\bibitem{L} A.~Leibman, \emph{Pointwise convergence of ergodic averages for polynomial
sequences of translations on a nilmanifold}, Ergodic Theory Dynam. Systems
\textbf{25} (2005), 201--213.

\bibitem{sarnak-lectures} P.~Sarnak, \emph{Three lectures on the M\"obius Function randomness and dynamics}, 2010, see
http://publications.ias.edu/sarnak/paper/506.


\bibitem{TT-blog} T.~Tao, \emph{The Chowla conjecture and the Sarnak conjecture}, see https://terrytao.wordpress.com/2012/10/14/the-chowla-conjecture-and-the-sarnak-conjecture/.

\bibitem{Z} T.~Ziegler, \emph{A soft proof of orthogonality of M\"obius to nilflows}, see http://tx.technion.ac.il/~tamarzr/soft-green-tao.pdf.


\end{thebibliography}
\end{document}